\input amstex
\input amsppt.sty
 \magnification=\magstep1
 \hsize=30truecc
 \vsize=22.2truecm
 \baselineskip=16truept
 \nologo
 \TagsOnRight
 \def\N{\Bbb N}
 \def\Z{\Bbb Z}
 
 \def\Q{\Bbb Q}
 
 \def\C{\Bbb C}
 \def\l{\left}
 \def\r{\right}
 \def\bg{\bigg}
 \def\({\bg(}
 \def\[{\bg[}
 \def\){\bg)}
 \def\]{\bg]}
 \def\t{\text}
 \def\f{\frac}
 \def\colon{{:}\;}
 \def\ch{\roman{ch}}
 
 \def\per{\roman{per}}
 \def\em{\emptyset}
 \def\se {\subseteq}

 \def\bi{\binom}
 
 \def\cs{\ldots}
 \def\ls{\leqslant}
 \def\gs{\geqslant}
 \def\al{\alpha}
 \def\ve{\varepsilon}
 \def\da{\delta}

 \def\bi{\binom}
 
 \def\Proof{\noindent{\it Proof}}
 
 \def\Remark{\medskip\noindent{\it  Remark}}
 \def\Ack{\medskip\noindent {\bf Acknowledgments}}

 \hbox{J. Combin. Theory Ser. A 103(2003), no.2, 288--301.}
 \medskip
  \topmatter
  \title On Snevily's Conjecture and Restricted Sumsets\endtitle
  \author Zhi-Wei Sun \endauthor
 \affil Department of Mathematics
     \\Nanjing University
     \\Nanjing 210093
    \\The People's Republic of China
    \\{\it E-mail:} {\tt zwsun\@nju.edu.cn}
    \\Homepage: {\tt http://pweb.nju.edu.cn/zwsun}
 \endaffil
 \abstract Let $G$ be an additive abelian group whose finite
 subgroups are all cyclic. Let $A_1,\cs,A_n$ $(n>1)$ be finite
 subsets of $G$ with cardinality $k>0$, and let $b_1,\cs,b_n$
 be pairwise distinct elements of $G$ with odd order. We show that for
 every positive integer $m\ls(k-1)/(n-1)$ there are
 more than $(k-1)n-(m+1)\bi n2$ sets $\{a_1,\cs,a_n\}$
 such that $a_1\in A_1,\cs, a_n\in A_n$, and both $a_i\not=a_j$ and
 $ma_i+b_i\not=ma_j+b_j$ (or both $ma_i\not=ma_j$ and $a_i+b_i\not=a_j+b_j$)
 for all $1\ls i<j\ls n$.
 This extends a recent result of Dasgupta, K\'arolyi, Serra and
 Szegedy on Snevily's conjecture. Actually  stronger
 results on sumsets with polynomial restrictions are obtained in
 this paper.
 \endabstract
 \thanks 2000 {\it Mathematics Subject Classification.}
Primary 11B75; Secondary 05A05, 11C20, 11P70, 12E05, 20D60.
\newline
\indent
Supported by
 the Teaching and Research Award Program for Outstanding Young Teachers
in Higher Education Institutions of MOE, and
 the National Natural
Science Foundation of the People's Republic of China.
 \endthanks
 \endtopmatter
  \document
 \heading {1. Introduction}\endheading

 In 1999 Snevily [Sn] posed the following conjecture.

\proclaim{Snevily's Conjecture} Let $G$ be an additive abelian
group with $|G|$ odd. Let $A$ and $B$ be subsets of $G$ with
cardinality $n>0$. Then there is a numbering $\{a_i\}_{i=1}^n$
of the elements of $A$ and a numbering $\{b_i\}_{i=1}^n$ of the elements
of $B$ such that $a_1+b_1,\cs,a_n+b_n$ are pairwise distinct.
\endproclaim

Using the polynomial method of Alon, Nathanson and Ruzsa (see,
e.g. [ANR], [A1] and [N]),
Alon [A2] proved that the above conjecture holds when $|G|$
is an odd prime.
In 2001 Dasgupta, K\'arolyi, Serra and
 Szegedy [DKSS] confirmed Snevily's conjecture for any cyclic group
 with odd order.

 In this paper we will show the following result in this direction.

 \proclaim{Theorem 1.1} Let $G$ be an additive abelian group whose finite
 subgroups are all cyclic. Let $A_1,\cs,A_n$ $(n>1)$ be finite
 subsets of $G$ with cardinality $k\gs n$, and let $b_1,\cs,b_n$
 be elements of $G$. Let $m$ be any positive integer
 not exceeding $(k-1)/(n-1)$.

 {\rm (i)} If $b_1,\cs,b_n$ are pairwise distinct,
 then there are at least $(k-1)n-m\bi n2+1$
 multisets $\{a_1,\cs,a_n\}$ such that $a_i\in A_i$ for $i=1,\cs,n$
 and all the $ma_i+b_i$ are pairwise distinct.

 {\rm (ii)} The sets
  $$\{\{a_1,\cs,a_n\}\colon a_i\in A_i,\ a_i\not=a_j\
 \t{and}\ ma_i+b_i\not= ma_j+b_j\ \t{if}\ i\not=j\}\tag1.1$$
 and $$\{\{a_1,\cs,a_n\}\colon a_i\in A_i,\ ma_i\not=ma_j\
 \t{and}\ a_i+b_i\not= a_j+b_j\ \t{if}\ i\not=j\}\tag1.2$$
 have more than $(k-1)n-(m+1)\bi n2\gs(m-1)\bi n2$ elements,
 provided that $b_1,\cs,b_n$ are pairwise distinct and of odd order,
 or they have finite order and
 $n!$ cannot be written in the form $\sum_{p\in P}px_p$ where all the $x_p$
 are nonnegative integers and $P$ is the set of primes dividing
 one of the orders of $b_1,\cs,b_n$.
 \endproclaim

 \Remark\ 1.1. When $G$ is a cyclic group with $|G|$ being odd or
 a prime power,
 our Theorem 1.1 (ii) in the case $k=n$ and $m=1$,
 yields Theorems 1 and 2 of [DKSS] respectively.
 In our opinion, the condition that all finite subgroups of $G$
 are cyclic might be omitted from Theorem 1.1.
\medskip

 We will deduce Theorem 1.1 from our stronger
results on sumsets with polynomial restrictions.
(As for sumsets of subsets of $\Z$ with linear restrictions,
the reader may consult [Su2].)

 Let $F$ be a field. We use $\ch(F)$ to denote
 the additive order of the multiplicative
identity of $F$ and call it the {\it characteristic} of $F$.
(When $\ch(F)=\infty$, some mathematicians regard the
characteristic of $F$ as zero.)
There are several recent results ([DH], [ANR], [HS], [LS], [PS])
concerning various restricted sumsets of $A_1,\cs,A_n\se F$.
For example, Corollary 1 of [HS] in the case $m=1$ can be
stated as follows:

{\it  Let $k\gs n\gs1$ be integers,
and $F$ be a field with
$\ch(F)$ greater than $n$ and $(k-n)n$. Let $A_1,\cs, A_n$ be subsets of $F$
with cardinality $k$, and $b_1,\cs,b_n$ be elements of $F$.
Then the sumset
$$\{a_1+\cs+a_n\colon a_i\in A_i,\ a_i\not=a_j\ \t{and}\
a_i+b_i\not=a_j+b_j\ \t{if}\ i\not=j\}$$
has more than $(k-n)n$ elements.}

When $F$ is a finite field of order $p^{\al}$
(where $p$ is a prime and $\al$ is a positive integer),
the additive group of $F$ is isomorphic to  the direct sum of
$\al$ copies of the additive cyclic group $\Z/p\Z$,
and the above result in the case $k=n$ was also found
by Dasgupta et al. ([DKSS]) who followed Alon's approach in [A2].

Let $R$ be any commutative ring with identity.
For $P(x_1,\cs,x_n)\in R[x_1,\cs,x_n]$, we simply
write $[x_1^{i_1}\cdots x_n^{i_n}]P(x_1,\cs,x_n)$ for the coefficient
of the monomial $x_1^{i_1}\cdots x_n^{i_n}$ in $P(x_1,\cs,x_n)$.
For a matrix $A=(a_{ij})_{1\ls i,j\ls n}$ with
entries in $R$, we use
 $\det(A)$ or $|a_{ij}|_{1\ls i,j\ls n}$ to denote the determinant of $A$,
and define the {\it permanent} of $A$ by
$$\roman{per}(A)=\sum_{\sigma\in S_n}a_{1,\sigma(1)}\cdots
a_{n,\sigma(n)}\tag1.3$$
where $S_n$ is the symmetric group of all the permutations on
$\{1,\cs,n\}$.

By studying certain
coefficients of some related polynomials
in the next section, we are able to
 prove the following main theorems in Section 3.

 \proclaim{Theorem 1.2} Let $k,m,n$ be positive integers with
$k>m(n-1)$, and let $A_1,\cs,A_n$ be subsets of a field $F$ with
cardinality $k$. Let $K=(k-1)n-m\bi n2$
and $P_1(x),\cs,P_n(x)\in F[x]$ have degree $m$.

{\rm (i)} If $\ch(F)>K$ and all the $b_i=[x^m]P_i(x)\ (i=1,\cs,n)$ are pairwise
distinct,  then $|S|\gs K+1$ where
$$S=\bg\{\sum_{i=1}^na_i\colon a_1\in A_1,\cs,a_n\in A_n,\
 \t{and}\ P_i(a_i)\not=P_j(a_j)\ \t{if}\ i\not=j\bg\}.\tag1.4$$

{\rm (ii)} If $\ch(F)>K-\bi n2$ and the permanent of
$B=(b_j^{i-1})_{1\ls i,j\ls n}$ does not vanish, then $|T|\gs K-\bi n2+1$
where
$$T=\bg\{\sum_{i=1}^na_i\colon a_i\in A_i,\ a_i\not=a_j\
 \t{and}\ P_i(a_i)\not=P_j(a_j)\ \t{if}\ i\not=j\bg\}.\tag1.5$$

{\rm (iii)} We have $\per(B)\not=0$, if $F$ is the complex field $\C$,
$b_1,\cs,b_n$ are $q$th roots of unity,
 and $n!$ does not belong to the set
$$ D(q)=\bg\{\sum_{p\mid q}px_p\colon x_p\in\{0,1,2,\cs\}\ \t{for any prime
divisor}\ p\ \t{of}\ q\bg\}.\tag1.6$$
\endproclaim

\Remark\ 1.2. Let $b_1,\cs,b_n$ be pairwise distinct elements of a field $F$,
and let $B$ be the Vandermonde matrix $(b_j^{i-1})_{1\ls i,j\ls n}$.
If $\ch(F)=2$,
then $$\per(B)=\det(B)=\prod_{1\ls i<j\ls n}(b_j-b_i)\not=0$$
as observed by Dasgupta et al. [DKSS]. If
$F=\C$, $\per(B)=0$ and all the $b_i$ are $q$th roots of unity, then
$\prod_{1\ls i<j\ls n}(1-b_i/b_j)=2\omega$
for some algebraic integer $\omega\in E=\Q(e^{2\pi i/q})$
and hence $q$ must be even (otherwise the norms of those $1-b_i/b_j$
($1\ls i<j\ls n$) with respect to the field extension $E/\Q$ would be
odd); this fact is also due to Dasgupta et al. [DKSS].
\medskip

Concerning the sumset $S$ given by (1.4), there is another result due to Liu and
Sun [LS]: {\it Let $F$ be a field and $P_1(x),\cs,P_n(x)\in F[x]$
be monic and of degree $m>0$. Let $A_1,\cs,A_n$ be finite subsets
of $F$ with $k=|A_n|>m(n-1)$ and $|A_{i+1}|-|A_i|\in\{0,1\}$ for
all $i=1,\cs,n-1$. If $L=(k-1)n-(m+1)\bi n2<\ch(F)$,
then $|S|\gs L+1$.}

\proclaim{Theorem 1.3} Let $A_1,\cs,A_n$ be finite subsets of
 a field $F$ with $0<k_1=|A_1|\ls\cs\ls k_n=|A_n|$, and let
 $P_1(x),\cs,P_n(x)\in F[x]$ be monic and of degree $m$ where
 $$m>k_n-k_1\ \t{and}\ k_n>m(n-1).\tag1.7$$

{\rm (i)} We have $L=\sum_{i=1}^n(k_i-1)-(m+1)\bi n2\gs0$.
If $\ch(F)>L!n!$, then $|T|\gs L+1$ where $T$ is as in $(1.5)$.

{\rm (ii)} When $k_1=\cs=k_n=k$ and $\ch(F)>L$, we have $|T|\gs L\gs(m-1)\bi
 n2$, and $|T|=L$ only if $k=n\gs \ch(F)>m=1$ or $\ch(F)=m=n=2<k=3$.
\endproclaim

As for Theorem 1.3 (ii), the following example
shows that $|T|=L$ may happen in the exceptional cases.

\medskip
\noindent {\it Example} 1.1. (i) Let $F$ be a field of
prime characteristic $p$. Let $b_1=0$ and $b_2=\cs=b_p=b\in
F\setminus\{0\}$.
Set $A_1=\cs=A_p=\{0,b,\cs,(p-1)b\}$.
Suppose that $a_1\in A_1,\cs,a_p\in A_p$ and $a_1,\cs,a_p$ are
pairwise distinct.
Then $a_1+b_1,\cs,a_p+b_p$ cannot be pairwise distinct.
In fact, if $a_1-b=a_i$ where $1<i\ls p$, then
$a_1+b_1=a_1=a_i+b_i$.

(ii) Let $F$ be a field of order 4 with identity 1.
Let $a\in F\setminus\{0,1\}$. Since $a^3=1$ and $a\ne1$, we have $a^2+a+1=0$.
If $a_1,a_2$ are distinct elements of $\{0,1,a\}$
and $a_1^2+a_1\not=a_2^2+a_2+1$, then
$\{a_1,a_2\}\not=\{0,a\},\{1,a\}$ and hence
$a_1+a_2=0+1$.
\medskip

The following example indicates that the condition $k_n>m(n-1)$
in Theorem 1.3 cannot be replaced by $k_n\gs
m(n-1)$.

\medskip
\noindent{\it Example} 1.2. Let $m$ and $n$ be positive integers.
Let $F$ be a finite field with $|F|=p^{\varphi(m)}>m(n-1)$,
where $p$ is a prime not dividing $m$ and $\varphi$
is Euler's totient function. As $m\mid p^{\varphi(m)}-1$,
the cyclic group $F^*=F\setminus\{0\}$ contains an element
$\gamma$ of order $m$. Since $|\{a^m\colon a\in F^*\}|=|F^*|/m\gs n-1$,
there are $c_1,\cs,c_{n-1}\in F^*$ such that $c_1^m,\cs,c_{n-1}^m$
are pairwise distinct. Clearly the set
$A=\{c_i\gamma^j\colon 0<i<n,\ 0\ls j<m\}$
has cardinality $m(n-1)$. If $a_1,\cs,a_n\in A$, then
$\{a_1^m,\cs,a_n^m\}\se\{c_1^m,\cs,c_{n-1}^m\}$
and so $a_1^m,\cs,a_n^m$ cannot be pairwise distinct.
\medskip

Now we give one more theorem.

\proclaim{Theorem 1.4} Let $F$ be a field, $b_1,\cs,b_n\in F$ and $c_{ij}\in F$
for all $1\ls i<j\ls n$. Let $P_1(x),\cs,P_n(x)\in F[x]$ be monic and of degree
$m>0$, and let $A_1,\cs,A_n$ be subsets of $F$
with cardinality $k>m(n-1)$. Set $B=(b_j^{i-1})_{1\ls i,j\ls n}$ and
$$C=\bg\{\sum_{i=1}^na_i\colon a_i\in A_i,\ P_i(a_i)\not=P_j(a_j)\
\t{and}\ a_ib_i-a_jb_j\not=c_{ij}\ \t{if}\ i<j\bg\}.\tag1.8$$

{\rm (i)} If $\ch(F)>(k-1)n-(m+1)\bi n2$ and $\per(B)\not=0$,
then $|C|>(k-1)n-(m+1)\bi n2\gs(m-1)\bi n2$.

{\rm (ii)} If $\ch(F)=2,\ m=1,\ k=n+1$ and $b_1,\cs,b_n$ are pairwise
distinct, then we have $|C|\gs n+1$.

{\rm (iii)} Suppose that $F$ is the complex field and $b_1,\cs,b_n$
are $q$th roots of unity. If $n!\not\in D(q)$ where $D(q)$ is as in $(1.6)$,
or $q$ is odd and $b_1,\cs,b_n$ are pairwise distinct,
then $|C|\gs(k-1)n-(m+1)\bi n2+1$.
\endproclaim

\proclaim{Corollary 1.1 {\rm ([DKSS])}} Let $F$ be a field of characteristic $2$,
and let $A$ and $B=\{b_1,\cs,b_n\}$ be subsets of $F$ with cardinality
$n$. Then there is a numbering $\{a_i\}_{i=1}^n$ of the elements of
$A$ such that $a_1b_1,\cs,a_nb_n$ are pairwise distinct.
\endproclaim
\Proof. If $A=F$ then we may simply take $a_i=b_i$ because
$b_1^2,\cs,b_n^2$ are pairwise distinct.
If $a\in F\setminus A$, then we may apply Theorem 1.4 (ii) with
$A_1=\cs=A_n=A\cup\{a\}$. \qed

For an odd integer $n>0$, the multiplicative group of the finite
field $F$ with $|F|=2^{\varphi(n)}$ has a cyclic subgroup of order $n$.
This observation of Dasgupta et al. indicates that
 Corollary 1.1 implies the truth of Snevily's conjecture
for any cyclic group of odd order.

Later we will deduce Theorem 1.1 from Theorems 1.2 and 1.4.

\heading{2. Auxiliary Results}\endheading

For convenience
we set $\N=\{0,1,2,\cs\}$, $\Z^+=\{1,2,3,\cs\}$
and $(x)_n=\prod_{0\ls j<n}(x-j)$ for $n\in\N$.
(An empty product is regarded as 1.)
For $\sigma\in S_n$ we let $\ve(\sigma)$ take $1$ or $-1$
 according to whether $\sigma\in S_n$ is even or odd.

\proclaim{Theorem 2.1} Let $R$ be a commutative ring with
identity, and let $a_{i,j}\in R$ for all $i,j=1,\cs,n$.
Let $k_1,\cs,k_n,m_1,\cs,m_n$ be nonnegative integers with
$M=\sum_{i=1}^n m_i+\da\bi n2\ls \sum_{i=1}^nk_i$ where $\da\in\{0,1\}$.
Then
$$\align&[x_1^{k_1}\cdots x_n^{k_n}]
|a_{i,j}x_j^{m_i}|_{1\ls i,j\ls n}\prod_{1\ls i<j\ls n}
(x_j-x_i)^{\da}\cdot\bg(\sum_{s=1}^nx_s\bg)^{\sum_{i=1}^n k_i-M}
\\&\quad\ \ =\cases\sum_{\sigma\in S_n,\
D_{\sigma}\se\N}\ve(\sigma)N_{\sigma}\prod_{i=1}^na_{i,\sigma(i)}&\t{if}\ \da=0,
\\\sum_{\sigma\in T_n}\ve(\sigma')N_{\sigma}
\prod_{i=1}^na_{i,\sigma(i)}&\t{if}\ \da=1,
\endcases\endalign$$
where
$$\align D_{\sigma}=&\{k_{\sigma(1)}-m_1,\cs,k_{\sigma(n)}-m_n\},
\\T_n=&\{\sigma\in S_n\colon D_{\sigma}\se\N\ \t{and}\ |D_{\sigma}|=n\},
\\N_{\sigma}=&\f{(k_1+\cs+k_n-M)!}{\prod_{i=1}^n\prod\Sb0\ls j<k_{\sigma(i)}-m_i
\\ j\not\in D_{\sigma}\ \t{if}\ \da=1\endSb(k_{\sigma(i)}-m_i-j)}\in\Z^+,
\endalign$$
and $\sigma'\ ($with $\sigma\in T_n)$ is the unique permutation in $S_n$ such that
$$0\ls k_{\sigma(\sigma'(1))}-m_{\sigma'(1)}<\cs<k_{\sigma(\sigma'(n))}-m_{\sigma'(n)}.$$
\endproclaim
\Proof. Write
$$\align P(x_1,\cs,x_n)=&|a_{i,j}x_j^{m_i}|_{1\ls i,j\ls n}
\prod_{1\ls i<j\ls n}(x_j-x_i)^{\da}
\\=&\sum\Sb i_1,\cs,i_n\in\N\\i_1+\cs+i_n=M\endSb c_{i_1,\cs,i_n}x_1^{i_1}\cdots x_n^{i_n}
\endalign$$
where $c_{i_1,\cs,i_n}\in R$. And let $c=[x_1^{k_1}\cdots
x_n^{k_n}]P(x_1,\cs,x_n)(x_1+\cs+x_n)^{K}$ where
$K=k_1+\cs+k_n-M$. Clearly
$$\align c=&[x_1^{k_1}\cs x_n^{k_n}]P(x_1,\cs,x_n)
 \sum\Sb j_1,\cs,j_n\in\N\\ j_1+\cs+j_n=K\endSb
 \f{K!}{j_1!\cdots j_n!}x_1^{j_1}\cdots x_n^{j_n}
 \\=&\sum\Sb 0\ls i_1\ls k_1,\cs,0\ls i_n\ls k_n\\i_1+\cs+i_n=M\endSb
 \f{K!}{(k_1-i_1)!\cdots(k_n-i_n)!}c_{i_1,\cs,i_n}.
 \endalign$$

 It is well known that
$$\prod_{1\ls i<j\ls n}(x_j-x_i)=|x_j^{i-1}|_{1\ls i,j\ls n}
=\sum_{\sigma\in S_n}\ve(\sigma)\prod_{j=1}^nx_j^{\sigma(j)-1}.$$
In the case $\da=1$, we have
$$\align P(x_1,\cs,x_n)=&\sum_{\sigma\in S_n}\ve(\sigma)\prod_{j=1}^nx_j^{\sigma(j)-1}
\times\sum_{\tau\in S_n}\ve(\tau)\prod_{j=1}^na_{\tau(j),j}x_j^{m_{\tau(j)}}
\\=&\sum_{\sigma\in S_n}\ve(\sigma)\sum_{\tau\in S_n}\ve(\tau)
\prod_{j=1}^na_{\tau(j),j}x_j^{\sigma(j)-1+m_{\tau(j)}}
\endalign$$
and hence
$$\align c=&\sum\Sb\sigma,\tau\in S_n\\\sigma(j)-1+m_{\tau(j)}\ls k_j
\\\t{for}\ j=1,\cs,n\endSb\ve(\sigma)\ve(\tau)\f{K!}
{\prod_{j=1}^n(k_j-(\sigma(j)-1+m_{\tau(j)}))!}\prod_{j=1}^na_{\tau(j),j}
\\=&\sum\Sb\sigma,\tau\in S_n\\\sigma(j)-1+m_{\tau(j)}\ls k_j
\\\t{for}\ j=1,\cs,n\endSb\ve(\sigma)\ve(\tau)
\f{K!\prod_{j=1}^n(k_j-m_{\tau(j)})_{\sigma(j)-1}}
{\prod_{j=1}^n(k_j-m_{\tau(j)})!}\prod_{j=1}^na_{\tau(j),j}
\\=&\sum\Sb \tau\in S_n\\m_{\tau(j)}\ls k_j\\\t{for}\ j=1,\cs,n\endSb\ve(\tau)
\f{K!\sum_{\sigma\in S_n}\ve(\sigma)\prod_{j=1}^n(k_j-m_{\tau(j)})_{\sigma(j)-1}}
{\prod_{j=1}^n(k_j-m_{\tau(j)})!}\prod_{j=1}^na_{\tau(j),j}
\\=&\sum\Sb \tau\in S_n\\m_{\tau(j)}\ls k_j\\\t{for}\ j=1,\cs,n\endSb\ve(\tau)
\f{K!|(k_j-m_{\tau(j)})_{i-1}|_{1\ls i,j\ls n}}
{\prod_{j=1}^n(k_j-m_{\tau(j)})!}\prod_{j=1}^na_{\tau(j),j}.
\endalign$$
Similarly, if $\da=0$ then
$$\align c=&\sum\Sb \tau\in S_n\\m_{\tau(j)}\ls k_j\\\t{for}\ j=1,\cs,n\endSb\ve(\tau)
\f{K!}{\prod_{j=1}^n(k_j-m_{\tau(j)})!}\prod_{j=1}^na_{\tau(j),j}
\\=&\sum\Sb\sigma\in S_n\\D_{\sigma}\se\N\endSb\ve(\sigma)
\f{K!}{\prod_{i=1}^n(k_{\sigma(i)}-m_i)!}\prod_{i=1}^na_{i,\sigma(i)}
\endalign$$
as desired.

Since $x^r=(x)_r+\sum_{0\ls t<r}S(r,t)(x)_t$
for $r=0,1,\cs,n-1$ where $S(r,t)$ are Stirling numbers of
the second kind, for $\tau\in S_n$ we have
$$\align&|(k_j-m_{\tau(j)})_{i-1}|_{1\ls i,j\ls n}
=|(k_j-m_{\tau(j)})^{i-1}|_{1\ls i,j\ls n}
\\=&\prod_{1\ls s<t\ls n}\l(k_t-m_{\tau(t)}-(k_s-m_{\tau(s)})\r)
\\=&(-1)^{|\{1\ls s<t\ls n\colon
\tau(s)>\tau(t)\}|}\prod\Sb 1\ls s,t\ls n\\\tau(s)<\tau(t)\endSb
\l(k_t-m_{\tau(t)}-(k_s-m_{\tau(s)})\r)
\\=&\ve(\tau)\prod_{1\ls i<j\ls
n}\l(k_{\tau^{-1}(j)}-m_j-(k_{\tau^{-1}(i)}-m_i)\r).
\endalign$$
Therefore, if $\da=1$ then
$$\align c=&\sum\Sb \sigma\in S_n\\D_{\sigma}\se\N\endSb
\f{K!\prod_{1\ls i<j\ls n}(k_{\sigma(j)}-m_j-(k_{\sigma(i)}-m_i))}
{\prod_{i=1}^n(k_{\sigma(i)}-m_i)!}\prod_{i=1}^na_{i,\sigma(i)}
\\=&\sum_{\sigma\in T_n}\f{K!\ve(\sigma')\prod_{u,v\in D_{\sigma},\ u<v}(v-u)}
{\prod_{i=1}^n\prod_{0\ls j<k_{\sigma(i)}-m_i}(k_{\sigma(i)}-m_i-j)}
\prod_{i=1}^na_{i,\sigma(i)}
\\=&\sum_{\sigma\in T_n}\ve(\sigma')N_{\sigma}\prod_{i=1}^na_{i,\sigma(i)}.
\endalign$$
By the above we also have $N_{\sigma}\in\Z^+$ for all those
$\sigma\in T_n$.

The proof of Theorem 2.1 is now complete. \qed

\proclaim{Corollary 2.1} Let $R$ be a commutative ring with
identity, and let $A=(a_{ij})_{1\ls i,j\ls n}$ be a matrix with
all the $a_{ij}$ in $R$. Let $k,m_1,\cs,m_n$ be nonnegative
integers with $m_1\ls\cs\ls m_n\ls k$.

{\rm (i)} We have
$$\aligned&[x_1^k\cdots x_n^k]|a_{ij}x_j^{m_i}|_{1\ls i,j\ls
n}(x_1+\cs+x_n)^{kn-\sum_{i=1}^n m_i}
\\&=\f{(kn-\sum_{i=1}^n m_i)!}{\prod_{i=1}^n(k-m_i)!}\det(A).
\endaligned\tag2.1$$

{\rm (ii)} If $m_1<\cs<m_n$ then
$$\align &[x_1^k\cdots x_n^k]|a_{ij}x_j^{m_i}|_{1\ls i,j\ls n}
\prod_{1\ls i<j\ls n}(x_j-x_i)\cdot\bg(\sum_{s=1}^nx_s\bg)^{kn-\bi n2-\sum_{i=1}^n m_i}
\\&\quad\ =(-1)^{\bi n2}\f{(kn-\bi n2-\sum_{i=1}^n m_i)!}
{\prod_{i=1}^n\prod\Sb m_i<j\ls k\\j\not= m_{i+1},\cs,m_n\endSb(j-m_i)}\per(A).
\endalign$$
\endproclaim

\Proof. Let $\da\in\{0,1\}$, and suppose that $m_1<\cs<m_n$ if
$\da=1$. Then $K=kn-\sum_{i=1}^nm_i-\da\bi n2\gs 0$.
Set
$$c=[x_1^k\cdots x_n^k]|a_{ij}x_j^{m_i}|_{1\ls i,j\ls
n}\prod_{1\ls i<j\ls n}(x_j-x_i)^{\da}\cdot(x_1+\cs+x_n)^K.$$
By Theorem 2.1, if $\da=0$ then
$c=(K!/\prod_{i=1}^n(k-m_i)!)\det(A)$.
In the case $\da=1$ we should have
$$\align c=&\f{K!}{\prod_{i=1}^n
\prod\Sb 0\ls r<k-m_i\\k-r\not=m_{i+1},\cs,m_n\endSb(k-m_i-r)}
\sum_{\sigma\in S_n}\ve(\sigma')\prod_{i=1}^n a_{i,\sigma(i)}
\\=&\f{K!}{\prod_{i=1}^n\prod\Sb m_i<j\ls k\\j\not=m_{i+1},\cs,m_n\endSb(j-m_i)}
(-1)^{\bi n2}\per(A).
\endalign$$
This completes the proof. \qed

\proclaim{Corollary 2.2} Let $k_1,\cs,k_n,m_1,\cs,m_n$ be nonnegative integers
with
$$k_1-m_1>\cs>k_n-m_n\gs0\tag2.2$$
and $$\min_{1\ls i<n}(m_{i+1}-m_i)\gs
\max_{1\ls i\ls n}k_i-\min_{1\ls i\ls n}k_i.\tag2.3$$
Let $A=(a_{ij})_{1\ls i,j\ls n}$ be a matrix with $a_{ij}\in\N$
and $\prod_{i=1}^ma_{ii}\not=0$. Put $L=\sum_{i=1}^n(k_i-m_i)-\bi n2$.
Then, for the coefficient $c$ of $x_1^{k_1}\cdots x_n^{k_n}$
in the polynomial
$$|a_{ij}x_j^{m_i}|_{1\ls i,j\ls n}\prod_{1\ls i<j\ls
n}(x_j-x_i)\cdot(x_1+\cs+x_n)^{L},$$
we have $0<(-1)^{n(n-1)/2}c\ls L!\per(A)$.
\endproclaim
\Proof. In view of (2.3), if $\sigma\in S_n$ then
$$k_{\sigma(1)}-m_1\gs k_{\sigma(2)}-m_2\gs\cs\gs
k_{\sigma(n)}-m_n.$$
Thus, for any $\sigma\in T_n$ we have $\ve(\sigma')=(-1)^{\bi n2}$
because $\sigma'(i)=n-i+1$ for $i=1,\cs,n$. Observe that $T_n$
contains the identity of $S_n$ and $N_{\sigma}\ls L!$ for all
$\sigma\in T_n$. Therefore
$(-1)^{\bi n2}c=\sum_{\sigma\in T_n}N_{\sigma}\prod_{i=1}^na_{i,\sigma(i)}$
is a positive integer not larger than $L!\per(A)$.
We are done. \qed

\proclaim{Theorem 2.2} Let $R$ be a commutative ring with
identity, and let $A=(a_{ij})_{1\ls i,j\ls n}$ be a matrix with all
the $a_{ij}$ in $R$. Let $k,l_1,\cs,l_n,m_1,\cs,m_n$ be
nonnegative integers with $K=kn-\sum_{i=1}^n(l_i+m_i)\gs0$. Then
$$\aligned&[x_1^k\cdots x_n^k]|a_{ij}x_j^{l_i}|_{1\ls i,j\ls
n}|x_j^{m_i}|_{1\ls i,j\ls n}(x_1+\cs+x_n)^K
\\=&[x_1^k\cdots x_n^k]|a_{ij}x_j^{m_i}|_{1\ls i,j\ls
n}|x_j^{l_i}|_{1\ls i,j\ls n}(x_1+\cs+x_n)^K.
\endaligned\tag2.4$$
\endproclaim
\Proof. Let $c_1$ and $c_2$ denote the left-hand side and the
right-hand side of (2.4) respectively.
Observe that
$$\align c_1=&[x_1^k\cdots x_n^k]\sum_{\sigma\in S_n}\ve(\sigma)
\prod_{j=1}^n(a_{\sigma(j),j}x_j^{l_{\sigma(j)}})
\sum_{\tau\in S_n}\ve(\tau)\prod_{j=1}^nx_j^{m_{\tau(j)}}\cdot
\bigg(\sum_{s=1}^nx_s\bigg)^K
\\=&\sum\Sb \sigma,\tau\in S_n\\l_{\sigma(j)}+m_{\tau(j)}\ls k
\\\t{for}\ j=1,\cs,n\endSb\ve(\sigma)\ve(\tau)
\bigg[\prod_{j=1}^nx_j^{k-l_{\sigma(j)}-m_{\tau(j)}}\bigg]
(x_1+\cs+x_n)^K\prod_{j=1}^na_{\sigma(j),j}
\\=&\sum\Sb \sigma,\tau\in S_n\\l_{\sigma(j)}+m_{\tau(j)}\ls k
\\\t{for}\ j=1,\cs,n\endSb\ve(\sigma)\ve(\tau)
\f{K!}{\prod_{j=1}^n(k-l_{\sigma(j)}-m_{\tau(j)})!}
\prod_{j=1}^na_{\sigma(j),j}.
\endalign$$
Similarly,
$$c_2=\sum\Sb \sigma,\tau\in S_n\\m_{\sigma(j)}+l_{\tau(j)}\ls k
\\\t{for}\ j=1,\cs,n\endSb\ve(\sigma)\ve(\tau)
\f{K!}{\prod_{j=1}^n(k-m_{\sigma(j)}-l_{\tau(j)})!}
\prod_{j=1}^na_{\sigma(j),j}.$$

If $l_i>k$ for some $i=1,\cs,n$, then both $c_1$ and $c_2$ vanish.
Now suppose that $k\gs\max_{1\ls i\ls n}l_i$. Then
$$\align c_1=&\sum_{\sigma,\tau\in S_n}\ve(\sigma)\ve(\tau)
\f{K!\prod_{j=1}^n(k-l_{\sigma(j)})_{m_{\tau(j)}}}
{\prod_{j=1}^n(k-l_{\sigma(j)})!}\prod_{j=1}^na_{\sigma(j),j}
\\=&\sum_{\sigma\in S_n}\ve(\sigma)
\f{K!\sum_{\tau\in S_n}\ve(\tau)\prod_{j=1}^n(k-l_{\sigma(j)})_{m_{\tau(j)}}}
{(k-l_1)!\cdots(k-l_n)!}\prod_{j=1}^na_{\sigma(j),j}
\endalign$$
and
$$\align c_2=&\sum_{\sigma,\tau\in S_n}\ve(\sigma)\ve(\tau)
\f{K!\prod_{j=1}^n(k-l_{\tau(j)})_{m_{\sigma(j)}}}
{\prod_{j=1}^n(k-l_{\tau(j)})!}\prod_{j=1}^na_{\sigma(j),j}
\\=&\sum_{\sigma\in S_n}\ve(\sigma)
\f{K!\sum_{\tau\in S_n}\ve(\tau)\prod_{j=1}^n(k-l_{\tau(j)})_{m_{\sigma(j)}}}
{(k-l_1)!\cdots(k-l_n)!}\prod_{j=1}^na_{\sigma(j),j}.
\endalign$$
Note that
$$\align &\sum_{\tau\in S_n}\ve(\tau)\prod_{j=1}^n(k-l_{\sigma(j)})_{m_{\tau(j)}}
\\=&\sum_{\tau\in S_n}\ve(\tau)\prod_{i=1}^n(k-l_{\sigma\tau^{-1}(i)})_{m_i}
=\ve(\sigma)|(k-l_j)_{m_i}|_{1\ls i,j\ls n}
\\=&\ve(\sigma)\sum_{\tau\in S_n}\ve(\tau\sigma^{-1})
\prod_{i=1}^n(k-l_{\tau\sigma^{-1}(i)})_{m_i}
=\sum_{\tau\in S_n}\ve(\tau)\prod_{j=1}^n(k-l_{\tau(j)})_{m_{\sigma(j)}}.
\endalign$$
So we have $c_1=c_2$. \qed

\heading{3. Proof of Theorems 1.1--1.4}\endheading

Theorem 1.2 of [A1] implies the following basic lemma.

\proclaim{Lemma 3.1 {\rm ([A1, Theorem 4.1; ANR, Theorem 2.1])}}
Let $A_1,\cs,A_n$ be finite subsets of a field $F$ with $k_i=|A_i|>0$
for $i=1,\cs,n$. Let $P(x_1,\cs,x_n)\in F[x_1,\cs,x_n]\setminus\{0\}$ and $\deg P
\ls\sum_{i=1}^n (k_i-1)$.
If $$[x_1^{k_1-1} \cdots x_n^{k_n-1}]P(x_1,\cs,x_n)
(x_1+\cs+x_n)^{\sum_{i=1}^n(k_i-1)-\deg P}\not=0,$$
then
$$|\{a_1+\cs+a_n\colon a_i\in A_i, \ P(a_1,\cs,a_n)\not=0\}|
\gs \sum_{i=1}^n(k_i-1)-\deg P+1.$$
\endproclaim

\noindent{\it Proof of Theorem 1.2}.
Applying Corollary 2.1 we find that
$$\align&[x_1^{k-1}\cdots x_n^{k-1}]
\prod_{1\ls i<j\ls n}(P_j(x_j)-P_i(x_i))\cdot(x_1+\cs+x_n)^K
\\=&[x_1^{k-1}\cdots x_n^{k-1}]
|b_j^{i-1}x_j^{(i-1)m}|_{1\ls i,j\ls n}(x_1+\cs+x_n)^K
\\=&\f{K!}{\prod_{i=1}^n(k-1-(i-1)m)!}|b_j^{i-1}|_{1\ls i,j\ls n}
\\=&\f{K!}{\prod_{r=0}^{n-1}(k-1-rm)!}\prod_{1\ls i<j\ls n}(b_j-b_i)
\endalign$$
and
$$\align&\[\prod_{s=1}^nx_s^{k-1}\]\prod_{1\ls i<j\ls n}(x_j-x_i)(P_j(x_j)-P_i(x_i))
\cdot(x_1+\cs+x_n)^{K-\bi n2}
\\=&\[\prod_{s=1}^nx_s^{k-1}\]|b_j^{i-1}x_j^{(i-1)m}|_{1\ls i,j\ls n}
\prod_{1\ls i<j\ls n}(x_j-x_i)\cdot(x_1+\cs+x_n)^{K-\bi n2}
\\=&(-1)^{\bi n2}\f{(K-\bi n2)!}
{\prod_{i=1}^n\prod\Sb(i-1)m<j\ls k-1\\j/m\not=i,\cs,n-1\endSb(j-(i-1)m)}
\per(B).
\endalign$$
In view of Lemma 3.1 we have parts (i) and (ii) of Theorem 1.2.

Now suppose that $F$ is the complex field and $b_1,\cs,b_n$ are
$q$th roots of unity. Then $\per(B)=\sum_{\sigma\in
S_n}b_{\sigma}$
where $b_{\sigma}=\prod_{i=1}^nb_i^{\sigma(i)-1}$ is a $q$th root
of unity.  For any integer $t$ relatively prime to
$q$, the cyclotomic field $\Q(e^{2\pi i/q})$
has an automorphism $\rho_t$ with $\rho_t(e^{2\pi i/q})
=e^{2\pi i t/q}$ and therefore
$$\sum_{\sigma\in S_n}b_{\sigma}^t=\sum_{\sigma\in S_n}\rho_t(b_{\sigma})
=\rho_t(\per(B)).$$
If $\per(B)=0$, then $\sum_{\sigma\in S_n}b_{\sigma}^t=0$ for
all those $t\in\Z$ divisible by none of the
prime divisors of $q$, and thus $n!=|S_n|\in D(q)$ by Lemma 9 of
[Su1]. This proves part (iii) of Theorem 1.2.
\qed

\medskip
\noindent{\it Proof of Theorem 1.3}. (i) As $m>k_n-k_1\gs k_{i+1}-k_i$
for $i=1,\cs,n-1$, we have
$$k_1-1>k_2-1-m>\cs>k_n-1-(n-1)m\gs0.$$
Thus $\sum_{i=1}^n(k_i-1-(i-1)m)\gs\sum_{i=0}^{n-1}i$ and hence
$L\gs0$. Let
$$l=[x_1^{k_1-1}\cdots x_n^{k_n-1}]|x_j^{(i-1)m}|_{1\ls i,j\ls n}
\prod_{1\ls i<j\ls n}(x_j-x_i)\cdot(x_1+\cs+x_n)^L.$$
Then $0<(-1)^{n(n-1)/2}l\ls L!n!$ by Corollary 2.2. Observe that
$$[x_1^{k_1-1}\cdots x_n^{k_n-1}]\prod_{1\ls i<j\ls
n}(x_j-x_i)(P_j(x_j)-P_i(x_i))\cdot(x_1+\cs+x_n)^L$$
coincides with $le$ where $e$ is the identity of the field $F$.
If $\ch(F)>L!n!$, then $le\not=0$ and hence
$|T|\gs L+1$ by Lemma 3.1.

(ii) It is clear that
$$L\gs m(n-1)n-(m+1)\bi n2=(m-1)\bi n2=\f{(m-1)(n-1)}2n.$$
If $k_1=\cs=k_n=k$ and $n\gs \ch(F)>L$, then $k-1=m(n-1)$ and $(m-1)(n-1)<2$, thus
  $m=1$ and $k=n\gs \ch(F)$ (in this case $L=0$), or $m=n=2=\ch(F)$ and $k=3$
  (in this case $L=1$).

 In light of Theorem 1.2(ii) and the above, it suffices
 to deduce a contradiction under the conditions $\ch(F)=m=n=2$,
 $k_1=k_2=3$ and $T=\em$. Let $a,b,c$ be the three elements of
 $A_1$. If $d\in A_2\setminus A_1$, then $P_1(x)-P_2(d)=0$
 for $x=a,b,c$, which is absurd since $\deg P_1(x)=2$.
 Therefore $A_2\se A_1$ and hence $A_2=A_1=\{a,b,c\}$. As $T=\em$
 we have
 $P_1(a)=P_2(b)=P_1(c)=P_2(a)=P_1(b)$,
 thus $P_1(x)=P_2(a)$ for $x=a,b,c$, which also leads to a
 contradiction.

 The proof of Theorem 1.3 is now complete. \qed

 \medskip
 \noindent{\it Proof of Theorem 1.4}. Note that
 $L=(k-1)n-(m+1)\bi n2\gs(m-1)\bi n2$.
 In view of Theorem 2.2 and Corollary 2.1 (ii),
  $$\align&\bigg[\prod_{i=1}^nx_i^{k-1}\bigg]
  \prod_{1\ls i<j\ls n}(P_j(x_j)-P_i(x_i))(b_jx_j-b_ix_i+c_{ij})
  \cdot\bigg(\sum_{s=1}^nx_s\bigg)^L
  \\\ &=\bigg[\prod_{i=1}^nx_i^{k-1}\bigg]|x_j^{(i-1)m}|_{1\ls i,j\ls n}
  |b_j^{i-1}x_j^{i-1}|_{1\ls i,j\ls n}\bigg(\sum_{s=1}^nx_s\bigg)^L
  \\\ &=(-1)^{\bi n2}(N\per(B))
  \endalign$$
 where
 $$N=\f{L!}{\prod_{i=1}^n\prod\Sb (i-1)m<j\ls k-1\\j/m\not=i,\cs,n-1\endSb(j-(i-1)m)}.$$

 If $\ch(F)>L$ and $\per(B)\not=0$, then $N\per(B)\not=0$ and hence
 $|C|>L$ by Lemma 3.1. This proves part (i).

 When $\ch(F)=2,\ m=1,\ k=n+1$ and $b_1,\cs,b_n$ are pairwise distinct, we have
 $$\align N\per(B)=&\f{n!}{\prod_{i=1}^n
 \prod\Sb i-1<j\ls n\\j\not=i,\cs,n-1\endSb(j-(i-1))}\per(B)
  \\=&\per(B)=\det(B)=\prod_{1\ls i<j\ls n}(b_j-b_i)\not=0
 \endalign$$
 and hence $|C|\gs L+1=n+1$. So part (ii) also holds.

 Combining part (i) with Theorem 1.2 (iii) and Remark 1.2,
 we obtain part (iii) of Theorem 1.4. \qed

\medskip
\noindent{\it Proof of Theorem 1.1}. Let $H$ be the subgroup of
$G$ generated by the finite set $A_1\cup\cs\cup A_n\cup\{b_1,\cs,b_n\}$.
By the structure theorem for finitely generated abelian groups,
$H$ is isomorphic to the direct sum $\t{Tor}(H)\oplus \Z^r$ for some $r\in\N$,
where $\t{Tor}(H)=\{a\in H\colon \t{the order of}\ a \ \t{is
finite}\}$ is a finite subgroup of $G$ and hence cyclic.
Let $h=|\t{Tor}(H)|$ and
choose an even integer $h'>2$ so that $h\mid h'$ and
$\varphi(h')/2\gs r+1$. By Dirichlet's unit theorem (cf. [H, Theorem 100]),
the unit group $U_{h'}$ of the ring $\Z[e^{2\pi i/h'}]$
is isomorphic to $(\Z/h'\Z)\oplus\Z^{\varphi(h')/2-1}$. Thus we can
identify the additive group $H$ with a subgroup of
the multiplicative group $U_{h'}$.
So, without loss of generality,
we may simply let $G$ be
the multiplicative group
$\C^*=\C\setminus\{0\}$.
(There is an alternate way to embed $H$ into the group $\C^*$:
Take any minimal basis $c_0,\cs,c_r$ of $H$ where
$c_0$ is a torsion element of order $h$. Then map $c_0$
to $e^{2\pi i/h}$ and $c_1,\cs,c_r$ onto any set of $r$
algebraically independent elements of the uncountable unit circle
$\{e^{2\pi i\theta}\colon 0\ls \theta<1\}$.
Clearly the map can be extended in a unique way to an embedding of $H$
into the multiplicative group $\C^*$.)

(i) If $b_1,\cs,b_n\in\C^*$ are pairwise
distinct, then by Theorem 1.2 the sumset
$$\{a_1+\cs+a_n\colon a_1\in A_1,\cs, a_n\in A_n,\ \t{and}\
a_i^mb_i\not=a_j^mb_j\ \t{if}\ i\not=j\}$$
has at least $(k-1)n-m\bi n2+1$ elements.

(ii) Now let us work under the conditions of Theorem 1.1 (ii).
If $b_1,\cs,b_n\in\C^*$ are of finite order, then they are
$q$th roots of unity where $q$ is the least common multiple of the
multiplicative orders of $b_1,\cs,b_n$.
By Remark 1.2 and Theorems 1.2 and 1.4, both
$$\{a_1+\cs+a_n\colon a_i\in A_i,\ a_i\not=a_j\ \t{and}
\ a_i^mb_i\not=a_j^mb_j\ \t{if}\ i\not=j\}$$
and
$$\{a_1+\cs+a_n\colon a_i\in A_i,\ a_i^m\not=a_j^m\ \t{and}
\ a_ib_i\not=a_jb_j\ \t{if}\ i\not=j\}$$
have more than $(k-1)n-(m+1)\bi n2$ elements.

So far we have completed the proof of Theorem 1.1.
\qed

\Ack. The author thanks Prof. N. Alon for sending the paper
[DKSS] and the referees for their helpful comments.

\enddocument